\documentclass[a4paper,11pt]{article}
\usepackage{amsfonts}
\usepackage{amstext}
\usepackage{amsthm}
\usepackage{amsmath}
\usepackage{amssymb}
\usepackage{xcolor}
\usepackage{latexsym}
\usepackage{graphicx}
\usepackage[latin1]{inputenc}
\newcommand{\R}{\Bbb{R}}

\newtheorem{teor}{Theorem}[section]

\newtheorem{cor}{Corollary}[section]
\newtheorem{rem}{Remark}[section]

\newcommand{\n}{\noindent}

\begin{document}

\title{Eigenvalue estimates and maximum principle for Lane-Emden systems, and applications to poly-Laplacian equations
}

\author{\textbf{Sabri Bahrouni \footnote{\textit{E-mail addresses}:
sabri.bahrouni@fsm.rnu.tn (S. Bahrouni)}}\\ {\small\it Mathematics Department,
Faculty of Sciences, University of Monastir,}\\ {\small\it 5019 Monastir, Tunisia}\\
\textbf{Edir Júnior Ferreira Leite \footnote{\textit{E-mail addresses}:
edirleite@ufscar.br (E.J.F. Leite)}}\\ {\small\it Departamento de Matem\'{a}tica,
Universidade Federal de São Carlos,}\\ {\small\it 13565-905, São Carlos, SP, Brazil}\\
\textbf{Gustavo Ferron Madeira \footnote{\textit{E-mail addresses}:
gfmadeira@ufscar.br (G.F. Madeira)}}\\ {\small\it Departamento de Matem\'{a}tica,
Universidade Federal de São Carlos,}\\ {\small\it 13565-905, São Carlos, SP, Brazil}
}
\date{}{%{\it Preprint - December 12, 2009}}
 
\maketitle

\markboth{abstract}{abstract}
\addcontentsline{toc}{chapter}{abstract}

\hrule \vspace{0,2cm}

\n {\bf Abstract}

This paper deals with explicit upper and lower bounds for principal eigenvalues and the maximum principle associated to generalized Lane-Emden systems (GLE systems, for short). Regarding the bounds, we generalize the upper estimate of Berestycki, Nirenberg and Varadhan [Comm. Pure Appl. Math. (1994), 47-92] for the first eigenvalue of linear scalar problems on general domains to the case of strongly coupled GLE systems
with $m \geqslant 2$ equations on smooth domains. The explicit lower estimate we obtain is also used to derive a maximum principle to GLE systems relying  in terms of quantitative ingredients. Furthermore, as applications of the previous results, upper and lower estimates for the first eigenvalue of weighted poly-Laplacian eigenvalue problems with $L^p$ weights $(p>n)$ and Navier boundary condition are obtained. Moreover, a strong maximum principle depending on the domain and the weight function for scalar problems involving the poly-Laplacian operator is also established.
\vspace{0.5cm}
\hrule\vspace{0.2cm}

\section{Introduction}

We consider the following eigenvalue problem which is called generalized Lane-Emden system (GLE system, for short):

\begin{equation} \label{1.3}
\left\{
\begin{array}{rlllr}
-{\cal L} u &=& \Lambda \rho(x) S_\alpha(u) & {\rm in} & \Omega, \\

u&=&0 & {\rm on} & \partial \Omega,
\end{array}\right.
\end{equation}
where $\Omega\subset\R^n$ is a bounded open subset with $\partial \Omega \in C^{1,1}$, $u=(u_1,\ldots,u_m): \overline{\Omega} \to \mathbb{R}^m$, $m \geqslant 2$, ${\cal L} u =({\cal L}_1 u_1, \ldots, {\cal L}_m u_m)$, where for each $i=1, \ldots, m$
\[
{\cal L}_i = a^i_{\kappa l}(x) \partial_{\kappa l} + b^i_{l}(x) \partial_l + c_i(x),\ \ \kappa,l=1,\ldots,n
\] 
is a uniformly elliptic linear operator of second order in \textit{nondivergence} form in $\Omega$, $\Lambda = (\lambda_1,\ldots, \lambda_m) \in \R^m$, $\rho = (\rho_1, \ldots, \rho_m) \in L^p(\Omega; \R^m)$  positive function in $\Omega$, $p>n$ and $\Lambda \rho(x) S_\alpha(u) = (\lambda_1 \rho_1(x) S_{\alpha_1}(u), \ldots,$ $\lambda_m \rho_m(x) S_{\alpha_m}(u))$, where $S_\alpha(u)$ is given by

\[
(S_{\alpha_1}(u), \ldots, S_{\alpha_m}(u)) = (|u_2|^{\alpha_1-1}u_2, |u_3|^{\alpha_2-1}u_3, \ldots,|u_m|^{\alpha_{m-1}-1}u_m, |u_1|^{\alpha_m-1}u_1)
\]
with $\alpha = (\alpha_1, \ldots, \alpha_m) \in (0, \infty)^m$ and $\Pi \alpha := \Pi_{i = 1}^m \alpha_i = 1$.

We assume that the coefficients of each operator ${\cal L}_i$ satisfy $a^i_{\kappa l} \in C(\overline{\Omega})$ and there are constants $c_0, C_0>0$ and $b_0 \geqslant 0$ such that

\[
c_0 |\xi|^2 \leqslant a^i_{\kappa l}(x) \xi_\kappa \xi_l \leqslant C_0 |\xi|^2,\ \ \forall \xi \in \R^n,\ x \in \Omega
\]
and

\[
\left(\sum_{l = 1}^n |b^i_l(x)|^2 \right)^{1/2},\ |c_i(x)| \leqslant b_0,\ \ \forall x \in \Omega,\ \ \ \ i=1, \ldots, m. 
\]

We recall that the operator ${\cal L}_i$ satisfies the weak maximum principle ({\bf (WMP)}, for short) in $\Omega$ if, for any function $u \in W^{2,p}(\Omega)$ such that ${\cal L}_i u \leqslant 0$ in $\Omega$ and $u \geqslant 0$ on $\partial \Omega$, we have $u \geqslant 0$ in $\Omega$. If, in addition, either $u \equiv 0$ in $\Omega$ or $u > 0$ in $\Omega$, we say that the operator ${\cal L}_i$ satisfies the strong maximum principle ({\bf (SMP)}, for short) in $\Omega$. Equivalently, see \cite{BeNiVa} and \cite{LoMo}, the {\bf (SMP)} holds in $\Omega$ if, and only if, the first eigenvalue $\lambda_1({\cal L}_i,\Omega)$ of ${\cal L}_i$, for the Dirichlet problem, is positive, so \cite[Theorem 10.1]{CaLo} assures the validity of the {\bf (SMP)} for ${\cal L}_i$ in $\Omega$, $i = 1,\ldots, m$, under the following conditions:

\begin{equation}\label{hip5}
\Vert b^i\Vert^n_\infty\vert\Omega\vert\leqslant\left(c_0\sqrt{\Sigma}\right)^n\vert B_1\vert
\end{equation}
and
\begin{equation}\label{hip6}
c_0\Sigma\vert B_1\vert^{\frac{2}{n}}\vert\Omega\vert^{-\frac{2}{n}}-\Vert b^i\Vert_\infty\sqrt{\Sigma}\vert B_1\vert^{\frac{1}{n}}\vert\Omega\vert^{-\frac{1}{n}}-\inf_\Omega c_i>0,
\end{equation}
where $B_1$ is the unit ball of $\R^n$, $|\cdot |$ stands for the Lebesgue measure in $\R^n$,
\[
\Sigma:=\lambda_1(-\Delta, B_1)\ \ \ \ \text{and}\ \ \ \ \Vert b^i\Vert_\infty:=\sup_\Omega\left(\sum_{l=1}^{n}(b_l^i)^2\right)^\frac{1}{2}.
\]

For the GLE system \eqref{1.3}, we say that $\Lambda_0 = (\lambda_{01}, \ldots, \lambda_{0m}) \in \R^m$ is said to be an eigenvalue if the problem \eqref{1.3} admits a nontrivial strong solution $\varphi = (\varphi_1, \ldots, \varphi_m) \in W^{2,p}(\Omega; \R^m) \cap W_0^{1,p}(\Omega; \R^m)$ which is called an eigenfunction of \eqref{1.3} associated to $\Lambda_0$. The eigenvalue $\Lambda_0$ is said to be principal if the GLE system \eqref{1.3} admits an eigenfunction $\varphi$ such that all components $\varphi_i$ are positive in $\Omega$, $i=1,\ldots, m$. We say that $\Lambda_0$ is simple if for any other eigenfunction $\psi = (\psi_1, \ldots, \psi_m)$, there exists $\eta = (\eta_1, \ldots, \eta_m) \in \R^m\setminus\{0\}$ such that $\psi_i = \eta_i \varphi_i$ in $\Omega$ for all $i=1,\ldots, m$.

In \cite{Mo4}, Montenegro studied the existence of principal eigenvalues for \eqref{1.3} with $m=2$. In a constructive way and based on Leray-Schauder degree theory, Krasnoselskii's method and sub-supersolution technique to \eqref{1.3}, he established that the set of principal eigenvalues of \eqref{1.3} with $m=2$ determines a curve on first quadrant of $\R^2$ with some qualitative properties, such as simplicity, local isolation, continuity, asymptotic behaviour and min-max characterization.

For $m \geqslant 2$, as consequence of a nonlinear version of Krein-Rutman theorem, the authors in \cite{CdSM} showed that the set ${\bf \Lambda_1} = {\bf \Lambda_1}(-{\cal L}, \Omega, \rho)$ of principal eigenvalues of the GLE system \eqref{1.3} is given by the principal $(m-1)$-hypersurface ${\bf \Lambda_1}=\{\Lambda \in (0, \infty)^m :\, H(\Lambda) = \lambda_*\}$, where $\lambda_*=\lambda_*(\Omega)>0$ and $H: (0, \infty)^m \rightarrow \R$ is the function defined by $H(\Lambda)= H(\lambda_1, \ldots, \lambda_m) := \lambda_1 \lambda^{\alpha_1}_2 \ldots \lambda^{\alpha_1\ldots\alpha_{m-1}}_m.$ Thus,
\begin{equation}\label{hyp}
{\bf \Lambda_1}=\{\Lambda \in (0, \infty)^m :\, \lambda_1 \lambda^{\alpha_1}_2 \ldots \lambda^{\alpha_1\ldots\alpha_{m-1}}_m = \lambda_*\}
\end{equation}
(see Figure 1 for $m=3$). In Corollary 2.1 of \cite{CdSM}, the authors proved that for each vector $\sigma\in (0, \infty)^{m-1}$, there exists a unique number $\theta_*(\sigma)>0$ such that $ \theta_*(\sigma) (1, \sigma)\in {\bf \Lambda_1}$. Moreover, they showed that

\[
\theta_*(\sigma)=\left(\displaystyle{\frac{\displaystyle{\lambda_*}}{\prod_{j=1}^{m-1}\sigma_j^{\prod_{i=1}^j\alpha_i}}}\right)^{\displaystyle{
\frac{1}{\displaystyle{\sum_{j=1}^{m}\prod_{i=1}^j\alpha_i}}}}.
\]

\begin{figure}[ht]
\centering
\includegraphics[scale=0.8]{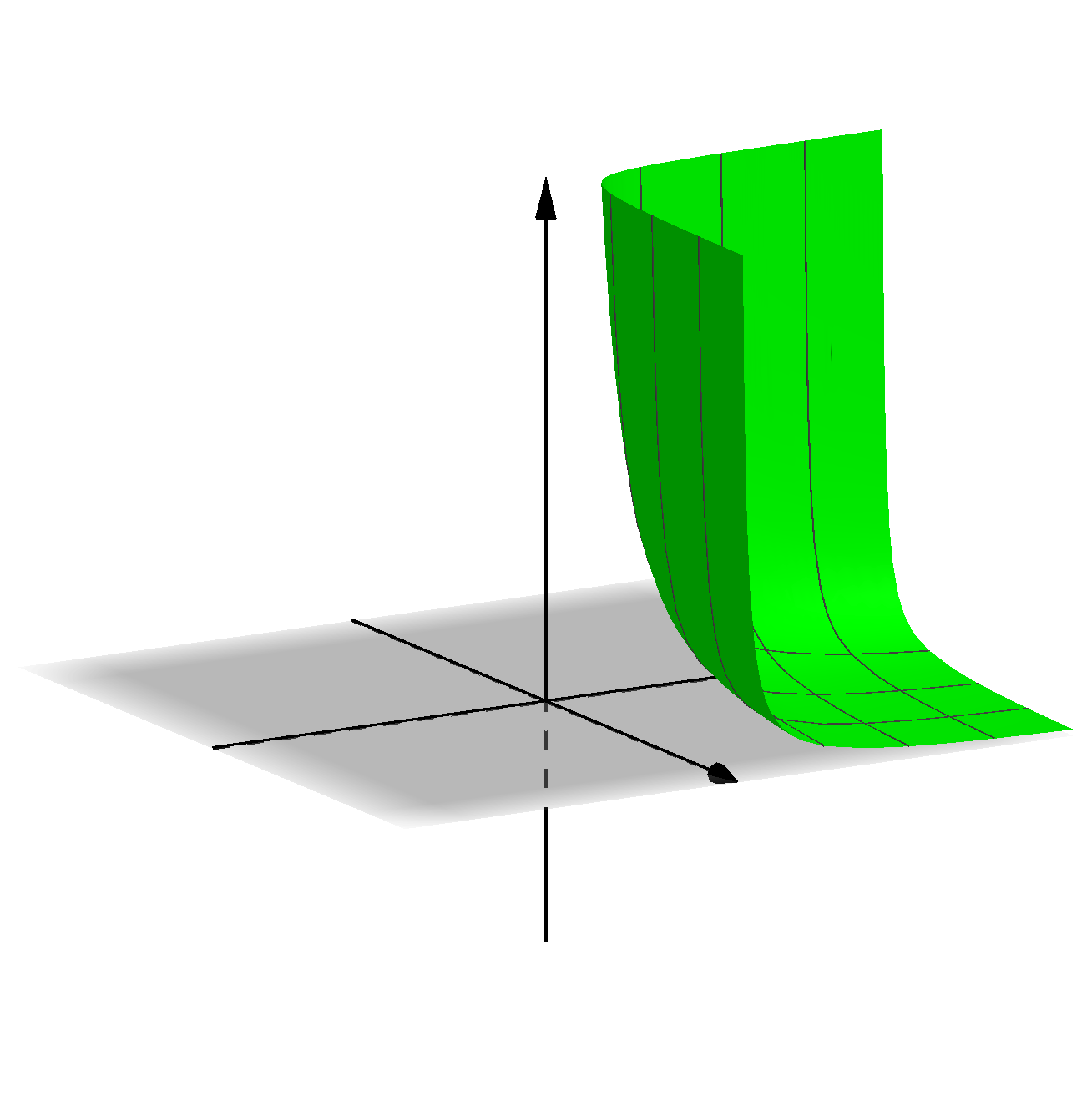}
\caption{Aspect of the principal $2$-hypersurface ${\bf \Lambda_1}=\left\{(\lambda_1,\lambda_2,\lambda_3) \in (0, \infty)^3 :\, \lambda_1 \lambda^{\alpha_1}_2\lambda^{\alpha_1\alpha_{2}}_3 = \lambda_*\right\}$.}
\end{figure}

In \cite{LM1}, Leite and Montenegro showed that the set ${\bf \Lambda_1}$ also satisfies the important properties mentioned above and obtained in \cite{Mo4} for $m=2$. Moreover, they proved that set ${\bf \Lambda_1}$ satisfy monotonicity property with respect to $\Omega$ and $\rho$ (see Theorem 6.1 of \cite{LM1}).

We say that {\bf (WMP)} associated to \eqref{1.3} holds in $\Omega$ if for any strong supersolution $u \in W^{2,p}(\Omega; \R^m)$ of the problem \eqref{1.3} such that $u \geqslant 0$ in $\partial\Omega$, we obtain $u \geqslant 0$ in $\Omega$. If furthermore either $u\equiv 0$ in $\Omega$ or, at least, one component of $u$ is positive in $\Omega$, then we say that {\bf (SMP)} associated to \eqref{1.3} holds in $\Omega$. If $\Lambda > 0$, then {\bf (SMP)} in $\Omega$ means that either $u\equiv 0$ in $\Omega$ or $u > 0$ in $\Omega$.

In Theorem 1.1 of \cite{LM1}, the authors characterize completely the set of $\Lambda \in \R^m$ such that {\bf (WMP)} and {\bf (SMP)} associated to \eqref{1.3} hold in $\Omega$. Precisely:

Let ${\cal R}_1$ be the open region in the first $m$-octant below ${\bf \Lambda_1}$, that is, ${\cal R}_1 := \{t \Lambda: 0 < t < 1\ {\rm and}\ \Lambda \in {\bf \Lambda_1}\}$. Then
\begin{equation}\label{carac}
\Lambda \in \overline{{\cal R}_1} \setminus {\bf \Lambda_1} \Leftrightarrow\ {\bf (WMP)}\text{ holds in}\ \Omega \Leftrightarrow\ {\bf (SMP)}\text{ corresponding to \eqref{1.3} holds in}\ \Omega.
\end{equation}

Moreover, using the ABP estimate (see \cite{GiTr}) they proved a lower estimate associated to (\ref{1.3}), and characterized the measure of $\Omega$ (depending on the diameter of $\Omega$ and the ABP constant) for which the maximum principles corresponding to (\ref{1.3}) hold in $\Omega$.

\subsection{Main results}

Regarding the {\bf (WMP)} and {\bf (SMP)} associated to nonlinear strongly coupled systems (\ref{1.3}), we are interested in establish explicitly the validity of such principles in terms of quantitative estimates involving the ingredients in (\ref{1.3}), such as the coefficients of $\mathcal{L}$, the weight $\rho\in L^p(\Omega; \R^m)$ and the Lebesgue measure $\vert\Omega\vert$. For this purpose, we shall obtain explicit lower estimates via Faber-Krahn's inequality. This famous inequality was obtained by Faber in \cite{faber} and Krahn in \cite{krahn}.

We define
\[
D:=\max\left\{\Vert\rho_i\Vert_{L^\infty(\Omega)}:i = 1,\ldots,m\right\},\ M:=\max\left\{\Vert b^i\Vert_{\infty}:i = 1,\ldots,m\right\}\ \ \text{and}\ \ \beta_0:={\displaystyle \sum_{i = 1}^m }\vert\inf_{\Omega}c_i\vert.
\]

\begin{teor}\label{teo1}($L^\infty$-weights) Let $\rho \in L^\infty(\Omega;\R^m)$. Suppose that $\Omega$ satisfies 
\begin{equation}\label{4.21}
M^n\vert\Omega\vert\leqslant\left(\frac{c_0\sqrt{\Sigma}}{\sqrt{m}}\right)^n\vert B_1\vert
\end{equation}
and \eqref{hip6} for $i = 1,\ldots,m$ and let ${\bf \Lambda_1}$ be the principal $(m-1)$-hypersurface (\ref{hyp}). Assume also that $\alpha_i=1$ for all $i = 1,\ldots,m$ (i.e., linear case). Then for $\Lambda_0=(\lambda_{01},\ldots,\lambda_{0m})\in {\bf \Lambda_1}$, we obtain
\begin{equation}\label{4.1}
{\displaystyle \sum_{i = 1}^m }\lambda_{0i}\geqslant\frac{1}{D}\left(c_0\Sigma\vert B_1\vert^{\frac{2}{n}}\vert\Omega\vert^{-\frac{2}{n}}  - M\sqrt{m}\sqrt{\Sigma} \vert B_1\vert^{\frac{1}{n}}\vert\Omega\vert^{-\frac{1}{n}} -\beta_0\right).
\end{equation}

Moreover, if in addition $\Omega$ satisfies 
\begin{equation}\label{01}
c_0\Sigma\vert B_1\vert^{\frac{2}{n}}\vert\Omega\vert^{-\frac{2}{n}}  - M\sqrt{m}\sqrt{\Sigma} \vert B_1\vert^{\frac{1}{n}}\vert\Omega\vert^{-\frac{1}{n}} -\beta_0> mD, 
\end{equation}
then $\lambda_*$ in (\ref{hyp}) satisfies
\begin{equation}\label{02}
\lambda_*\geqslant\frac{1}{D}\left(c_0\Sigma\vert B_1\vert^{\frac{2}{n}}\vert\Omega\vert^{-\frac{2}{n}}  - M\sqrt{m}\sqrt{\Sigma} \vert B_1\vert^{\frac{1}{n}}\vert\Omega\vert^{-\frac{1}{n}} -\beta_0\right).
\end{equation}
In particular,
\[
\lim_{\vert\Omega\vert\downarrow 0} \lambda_* =+\infty.
\]
\end{teor}

The above estimate was derived by refining the arguments employed in the proof of Theorem 5.1 of \cite{Lo} and in the proof of Theorem 4.1 of \cite{CauLo}. As a result of the preceding theorem and characterization (\ref{carac}), we obtain: 

\begin{cor} \label{EQ1} Under the hypotheses of Theorem \ref{teo1}, the following assertions are equivalent:

\begin{itemize}
\item[{\rm (i)}] $\Lambda \geqslant 0$;
\item[{\rm (ii)}] {\bf (WMP)} corresponding to \eqref{1.3} holds in $\Omega$ provided that 
\begin{equation}\label{4.3}
c_0\Sigma\vert B_1\vert^{\frac{2}{n}}\vert\Omega\vert^{-\frac{2}{n}}  - M\sqrt{m}\sqrt{\Sigma} \vert B_1\vert^{\frac{1}{n}}\vert\Omega\vert^{-\frac{1}{n}} -\beta_0>D{\displaystyle \sum_{i = 1}^m }\lambda_{i}
\end{equation}
when $\Lambda \in (0, \infty)^m$; 
\item[{\rm (iii)}] {\bf (SMP)} corresponding to \eqref{1.3} holds in $\Omega$ provided that $\Omega$ satisfies (\ref{4.3}) when $\Lambda \in (0, \infty)^m$.
\end{itemize}
\end{cor}

The nonlinear case, that is, when there is $\alpha_i \neq 1$ for some $i\in\{1,\ldots,m\}$, the proof of lower estimate is more intricate. In this situation, we draw upon the techniques used in the proofs of Theorems 5.1 in \cite{Lo}, 10.1 in \cite{CaLo}, and 1.6 in \cite{LM}.

We define
\[
C:=\max\left\{1,\Vert\rho_i\Vert_{L^p(\Omega)}:i = 1,\ldots,m\right\},\ \ \tilde{c_0}:=\min\left\{\left(\frac{c_0}{2}\right)^{{\displaystyle \Pi_{j = 1}^{i} \alpha_j}}\frac{1}{(3C)^{(\Pi_{j = 1}^{i} \alpha_j)-1}}: i = 1,\ldots,m\right\},
\]
\[
c:=\max\left\{\Vert b^i\Vert_\infty+\vert\inf_\Omega c_i\vert:i=1,\ldots,m\right\},\ \ \ \beta:=\min\left\{\Pi_{j = 1}^{i} \alpha_j:i= 1,\ldots,m\right\},
\]
\[
\ \ \text{and}\ \ \overline{c_0}:=\max\left\{c^{\Pi_{j = 1}^{i} \alpha_j}: i= 1,\ldots,m\right\}.
\]

\begin{teor}($L^p$-weights)\label{LE} Let $\rho \in L^p(\Omega;\R^m)$, where $p>n$ and $C_n$ given in the inequalities (\ref{sup}), (\ref{3}) and (\ref{4}). Suppose that $\Omega$ satisfies \eqref{hip5} and \eqref{hip6} for $i = 1,\ldots,m$ and let ${\bf \Lambda_1}$ be the principal $(m-1)$-hypersurface (\ref{hyp}). Assume also that $\alpha_1=\max\{\alpha_i:i = 1,\ldots,m\}$,
\begin{equation}\label{hip0}
\Vert b^i\Vert_\infty\vert\Omega\vert^{\frac{p(\alpha_1-1)+2}{2p(\alpha_1+1)}}\leqslant \frac{c_0}{2},\ \ \ i=1,\ldots m .
\end{equation}
and
\begin{equation}\label{hip1}
\vert\Omega\vert\leqslant\min\left\{1,C_n^{\frac{n(\theta-1)}{2\theta}}\Sigma^\frac{n}{2}\vert B_1\vert\right\}.
\end{equation}
Then, 

\begin{itemize}

\item[{\rm (i)}] For $n=1$, $\alpha_1\geqslant 1$ and $\Omega$ satisfying
\begin{equation}\label{h1}
\frac{1}{C_n^{\beta-\theta\beta}}\tilde{c_0}\Sigma^{\theta\beta}\vert B_1\vert^{\frac{2\theta\beta}{n}}\vert\Omega\vert^{-\frac{2\theta\beta}{n}}-\overline{c_0}\vert\Omega\vert^{\frac{p\beta(\alpha_1-1)+2\beta}{2p(\alpha_1+1)}}> 2C,
\end{equation}
we have $\lambda_*$ in (\ref{hyp}) satisfies
\begin{equation}\label{8.108}
\lambda_*\geqslant\frac{1}{2C}\left(\frac{1}{C_n^{\beta-\theta\beta}}\tilde{c_0}\Sigma^{\theta\beta}\vert B_1\vert^{\frac{2\theta\beta}{n}}\vert\Omega\vert^{-\frac{2\theta\beta}{n}}-\overline{c_0}\vert\Omega\vert^{\frac{p\beta(\alpha_1-1)+2\beta}{2p(\alpha_1+1)}}\right),
\end{equation}
where $\theta=\frac{2(p-1)}{p(\alpha_1+1)}$;
\item[{\rm (ii)}] For $n=2$, $\alpha_1\geqslant 1$ and $\Omega$ satisfying
\begin{equation}\label{h2}
\frac{1}{C_2^{\beta-\theta\beta}}\tilde{c_0}\Sigma^{\theta\beta}\vert B_1\vert^{\theta\beta}\vert\Omega\vert^{-\frac{\beta}{\alpha_1+1}}-\overline{c_0}\vert\Omega\vert^{\frac{\alpha_1\beta}{2(\alpha_1+1)}}> 2C,
\end{equation}
we get $\lambda_*$ in (\ref{hyp}) satisfies
\[
\lambda_*\geqslant\frac{1}{2C}\left(\frac{1}{C_2^{\beta-\theta\beta}}\tilde{c_0}\Sigma^{\theta\beta}\vert B_1\vert^{\theta\beta}\vert\Omega\vert^{-\frac{\beta}{\alpha_1+1}}-\overline{c_0}\vert\Omega\vert^{\frac{\alpha_1\beta}{2(\alpha_1+1)}}\right),
\]
where $\theta=\frac{1}{2\alpha_1+1}$;

\item[{\rm (iii)}] For $n\geqslant 3$, $1\leqslant\alpha_1<\frac{n}{n-2}$ and $\Omega$ satisfying (\ref{h1}), we obtain $\lambda_*$ in (\ref{hyp}) satisfies (\ref{8.108}), where
\[
\frac{p-1}{p(\alpha_1+1)}=\frac{\theta}{2}+\frac{(1-\theta)(n-2)}{2n}.
\]
\end{itemize}
Moreover, in any dimension,
\[
\lim_{\vert\Omega\vert\downarrow 0} \lambda_* =+\infty.
\]
\end{teor}

As an interesting consequence of characterization (\ref{carac}) and Theorem \ref{LE}, we derive the following maximum principle associated to \eqref{1.3}:

\begin{cor}\label{cor2} Let $\rho \in L^p(\Omega;\R^m)$, where $p>n$ and $\alpha > 0$ be such that $\Pi \alpha = 1$. When $n\geqslant 3$, assume also that $\left(\frac{n-2}{n}\right)^{m-1} < \alpha_i < \frac{n}{n-2}$ for all $i=1,\ldots,m$. The following assertions are equivalent:

\begin{itemize}
\item[{\rm (i)}] $\Lambda \geqslant 0$;
\item[{\rm (ii)}] There is an explicit positive constant $\eta_0$ depending only on $n$, $c_0$, $C_0$, $b_0$, $\alpha$, $\Lambda$ and $\rho$ such that {\bf (WMP)} corresponding to \eqref{1.3} holds in $\Omega$ whenever $|\Omega| < \eta_0$ when $\Lambda \in (0, \infty)^m$;
\item[{\rm (iii)}] There is an explicit positive constant $\eta_0$ depending only on $n$, $c_0$, $C_0$, $b_0$, $\alpha$, $\Lambda$ and $\rho$ such that {\bf (SMP)} corresponding to \eqref{1.3} holds in $\Omega$ whenever $|\Omega| < \eta_0$ when $\Lambda \in (0, \infty)^m$.
\end{itemize}
\end{cor}

Theorem \ref{teo1} is a special case of Theorem \ref{LE}. However, in the linear case, the proof is more straightforward and does not require, for instance, the use of Sobolev embedding, or the interpolation inequality. Therefore, we will present them separately.

Another significant area of interest involves the upper bounds of principal eigenvalues. For $m = 1$, Pucci \cite{P} provided a sharp upper bound for the principal eigenvalues of the problem \eqref{1.3}, with $b_l^1 = c_1 = 0$, $\rho = 1$, and $\Omega$ being a ball. In a broader context, the upper bound for the first eigenvalue was derived in Lemma 1.1 of \cite{BeNiVa} on general domains. The case for $m = 2$ was subsequently addressed in \cite{LM2} on smooth domains.

Finally, by the monotonicity of principal eigenvalues of (\ref{1.3}) with respect to $\Omega$, that is, if $\hat{\Omega} \supset \tilde{\Omega}$, then $\lambda_*(\hat{\Omega}) \leqslant \lambda_*( \tilde{\Omega})$ (see Theorem 6.1 of \cite{LM1}), we deduce the following explicit upper bound for all $m\geqslant 2$:

\begin{teor}\label{upper}
Suppose $\Omega$ satisfies \eqref{hip5} and \eqref{hip6} for $i = 1,\ldots,m$ and $\alpha > 0$ be such that $\Pi \alpha = 1$ and $\rho \in L^\infty(\Omega;\R^m)$ with $\rho_i \geqslant \varepsilon_0$ in $\Omega$ for some constant $\varepsilon_0 > 0$, $i = 1,\ldots, m$. Let ${\bf \Lambda_1}$ be the principal $(m-1)$-hypersurface (\ref{hyp}). Suppose also that $0 \in \Omega$. Let $0<R < 1$ such that $\overline{B}_R \subset \Omega$ and 
\begin{equation}\label{hip2}
c_i(0)R^2 <16 nc_0
\end{equation}
for $i =1,\ldots, m$. Then $\lambda_*$ in (\ref{hyp}) satisfies
\[
\lambda_*\leqslant ER^{-2\displaystyle{\sum_{j=1}^{m}\prod_{i=1}^j\alpha_i}},
\]
where $B_R=\{x \in \R^n: |x| < R\}$ and $E=A_1 A^{\alpha_1}_2 \ldots A^{\alpha_1\ldots\alpha_{m-1}}_m$, with
\begin{equation}\label{hip3}
A_i=\left\{
\begin{array}{rlllr}
 \dfrac{2^{4\alpha_i-1}}{\varepsilon_0c_0^{2\alpha_i - 1}}  \left( 4n C_0+\frac{9}{4} b_0 \right)  \left(4n C_0+\frac{9}{4}b_0 + 2c_0 \right)^{2\alpha_i - 1}  & {\rm if} & \alpha_i \geqslant 1/2, \\

 \dfrac{2^{6\alpha_i-2}}{\varepsilon_0}(4n C_0 + \frac{9}{4}b_0)  & {\rm if} & \alpha_i<1/2.
\end{array}\right.
\end{equation}
\end{teor}

Note that Theorems \ref{teo1}, \ref{LE}, and \ref{upper} extend, for example, Theorem 5.1 of \cite{Lo}, Theorem 4.1 of \cite{CauLo}, Theorem 10.1 of \cite{CaLo}, Theorem 1.6 of \cite{LM}, Theorem III of \cite{P}, Lemma 1.1 of \cite{BeNiVa}, and Theorem 3.2 of \cite{LM2} to the GLE system \eqref{1.3}.

\subsection{Applications to poly-Laplacian equations}

It is known that the eigenfunctions associated with the first eigenvalue of the bi-Laplacian with Dirichlet boundary condition are not always one-signed functions (see Subsection 3.1.3 in \cite{GGS}). This is due to the lack of a maximum principle to that class of problems. Further, Faber-Krahn type inequalities are also more delicate in this case (see Subsection 1.2.3.4 in \cite{HP}). Nevertheless, the topic of obtaining lower and upper bounds on eigenvalues of the poly-Laplacian operator with arbitrary order and Dirichlet boundary condition have been widely investigated, see for instance \cite{A, CW, CW1, JLWX, LP, Pe}. 

On the other hand, the study of estimates for eigenvalues of higher order elliptic equations with Navier boundary conditions seems to be little addressed in the literature. Some contributions on this subject are also provided in this paper. To the best of our knowledge, the only work in this direction is due to Lin and Zhu \cite{LZ}. In that paper, the authors have proven, among many other results, a delicate upper estimate of nodal sets of eigenfunctions of the poly-Laplacian with Navier boundary condition on real analytic domains $\Omega$ via analytic estimates of Morrey-Nirenberg and Carleman estimates.

Consider the following Navier boundary value problem involving the poly-Laplacian operator

\begin{equation}\label{poly}
\left\{
\begin{array}{llll}
(-\Delta)^m v = \lambda\rho_1(x) v & {\rm in} \ \ \Omega,\\
v=(-\Delta)^j v=0,\ j = 1,\ldots, m-1 & {\rm on} \ \ \partial\Omega.
\end{array}
\right.
\end{equation}

\begin{teor} \label{teo2} Let $\rho_1 \in L^p(\Omega)$, with $p>n$. Then the first eigenvalue $\lambda_1((-\Delta)^m,\Omega,\rho_1)$ of the problem (\ref{poly}) is positive and admits eigenfunction $v\in W^{2m,p}(\Omega)$ satisfying $ (-\Delta)^j v>0$ in $\Omega$, $j = 0,\ldots, m-1$. Moreover, $\lambda_1((-\Delta)^m,\Omega,\rho_1)$ is simple, i.e., for any other solution $z$ of (\ref{poly}) with $\lambda=\lambda_1((-\Delta)^m,\Omega,\rho_1)$, there exists $\eta\in\mathbb{R}\setminus\{0\}$ such that $z = \eta v$ in $\Omega$.
\end{teor}

\begin{teor}\label{teo4}
Let $\Omega=B_s$ and $\rho_1 \in L^\infty(B_s)$, with $\rho_1 \geqslant \varepsilon_0$ in $\Omega$ for some constant $\varepsilon_0 > 0$ and $ms^2\Vert\rho_1\Vert_{L^\infty(B_s)}\leqslant\Sigma$. Then for all $0<R < 1$, with $R<s$, the following estimates hold

\begin{equation}\label{i4}
\frac{n^2}{4s^2D}\leqslant\frac{\Sigma}{s^2D}\leqslant\lambda_1((-\Delta)^m,B_s,\rho_1)\leqslant\left(\dfrac{64n(2n+1)}{\varepsilon_0R^2}\right)^m,
\end{equation}
where $D=\max\{1,\Vert\rho_1\Vert_{L^\infty(\Omega)}\}$.
\end{teor}
\begin{rem}
\textnormal{Theorem \ref{teo4} extends to the poly-Laplacian the following classical inequality for the Laplacian
\begin{equation}\label{cheeger}
\left(\frac{n}{2}\right)^2\leqslant \lambda_1(-\Delta,B_1)\leqslant\dfrac{64n(2n+1)}{R^2},
\end{equation}obtained from Cheeger's constant (see \cite{Cheeger}) and Lemma 1.1 of \cite{BeNiVa}. Indeed, \eqref{cheeger} is recoverd from \eqref{i4} in the case $\rho_1(x)=1,$ $s=1$ and $m=1.$}
\end{rem}

Finally, we derive the following strong maximum principle for the poly-Laplacian operator with Navier boundary condition.

\begin{teor}\label{cor3} Let $\rho_1 \in L^p(\Omega)$, where $p>n$ and $\lambda\geqslant 0$. Then if $v\in W^{2m,p}(\Omega)$ is a nontrivial solution of the problem
\begin{equation}\label{MP}
\left\{
\begin{array}{llll}
(-\Delta)^m v \geqslant \lambda\rho_1(x) v & {\rm in} \ \ \Omega,\\
v=(-\Delta)^j v=0,\ j = 1,\ldots, m-1 & {\rm on} \ \ \partial\Omega,
\end{array}
\right.
\end{equation}
there exists an explicit constant $\eta_0=\eta_0(n,\lambda,\rho_1) > 0$, depending only on $n$, $\lambda$ and $\rho_1$, such that $(-\Delta)^j v>0$ in $\Omega$, $j = 0,\ldots, m-1,$ provided that $\vert\Omega\vert<\eta_0$.
\end{teor}

\begin{rem}
\textnormal{(i) It is worth to mentioning that the strong maximum principle holds for \eqref{MP} on any bounded open domain $\Omega\subset\R^n$ with $\partial \Omega \in C^{1,1}$ provided that $0\leqslant\lambda<\lambda_1((-\Delta)^m,\Omega,\rho_1)$, see Theorem 1.1 of \cite{LM1}. The issue in Theorem \ref{cor3} is to provide a sufficient condition for the  validity of the strong maximum principle for any given $\lambda\geqslant 0$ in terms of a suitable control on the size of the domain.}
\medskip

\noindent
\textnormal{(ii) A similar result to Theorem \ref{cor3} has been obtained by De Figueiredo in \cite{djairo}, Proposition 1.1. To point out some differences, in \cite{djairo} the bounded domain $\Omega$ is arbitrary, but $\rho_1 \in L^\infty(\Omega)$ and the constant $\eta_0$ depends on the diameter of $\Omega$ and the ABP constant. Moreover, results like Proposition 1.1 of \cite{djairo} and Theorem \ref{cor3} above can be applied to the study of symmetry properties of solutions of higher order elliptic equations, see for instance \cite{CV, djairo}.} 

\end{rem}

The rest of this paper is organized as follows: Section 2 focuses on proving the lower bound in Theorem \ref{teo1} and establishing Corollary \ref{EQ1}. In Section 3, we present the lower bounds for the principal eigenvalues associated with GLE systems, where the weight functions belong to $L^p(\Omega; \R^m)$ with $p > n$, as stated in Theorem \ref{LE}. Section 4 provides the proof of the upper estimate on $\lambda_*$ in Theorem \ref{upper}. Theorem \ref{teo2} is proved in Section 5. Finally, Section 6 explores both lower and upper estimates for the principal eigenvalue of (\ref{poly}), leading to the proof of Theorem \ref{teo4}. Moreover, the strong maximum principle for the poly-Laplacian operator, as stated in Theorem \ref{cor3}, is also proved in Section 6.

\section{Proof of Theorem \ref{teo1} and Corollary \ref{EQ1}}

Let's start with the proof of Theorem \ref{teo1}. For this, we denote by $\varphi$ a positive eigenfunction corresponding to $\Lambda_0\in {\bf \Lambda_1}$, normalized so that 
\begin{equation}\label{4.2}
{\displaystyle \sum_{i = 1}^m } ||\varphi_i||^2_{L^2(\Omega)} = 1.
\end{equation}

Since $\varphi$ is a positive strong solution of the GLE system \eqref{1.3}, then multiplying by $\varphi_i$ the scalar problem
\[
-{\cal L}_i\varphi_i = \lambda_{0i}\rho_i(x)\varphi_{i+1}
\]
 and integrating in $\Omega$, we get

\[
\lambda_{0i}\int_{\Omega}\rho_i(x)\varphi_{i+1}\varphi_i\geqslant c_0\left\Vert\nabla\varphi_i\right\Vert_{L^2(\Omega)}^2-\int_\Omega\langle b^i,\nabla\varphi_i\rangle\varphi_i-\inf_\Omega c_i\int_\Omega\varphi_i^2
\]
(see, e.g., page 266 of \cite{Lo1}), for $i=1,\ldots,m$. Note that
\begin{eqnarray*}
\left\vert\int_\Omega\langle b^i,\nabla\varphi_i\rangle\varphi_i\right\vert &\leqslant & \Vert b^i\Vert_\infty\left\Vert\nabla\varphi_i\right\Vert_{L^{2}(\Omega)},
\end{eqnarray*}
\[
\int_{\Omega}\rho_i(x)\varphi_{i+1}\varphi_i\leqslant \Vert \rho_i\Vert_{L^\infty(\Omega)}
\]
and
\[
-\inf_\Omega c_i\int_\Omega\varphi_i^2\geqslant -\vert\inf_\Omega c_i\vert\Vert\varphi_i\Vert_{L^{2}(\Omega)}^2.
\]
Therefore,

\begin{equation}\label{ine10}
\lambda_{0i}\Vert \rho_i\Vert_{L^\infty(\Omega)}\geqslant c_0\left\Vert\nabla\varphi_i\right\Vert^2_{L^2(\Omega)}  - \Vert b^i\Vert_\infty\left\Vert\nabla\varphi_i\right\Vert_{L^2(\Omega)} -\vert\inf_\Omega c_i\vert\Vert\varphi_i\Vert_{L^{2}(\Omega)}^2.
\end{equation}

So, adding up inequalities (\ref{ine10}) for all $i=1,\ldots,m$ shows that

\[
D{\displaystyle \sum_{i = 1}^m }\lambda_{0i}\geqslant c_0{\displaystyle \sum_{i = 1}^m }\left\Vert\nabla\varphi_i\right\Vert^2_{L^2(\Omega)}  - M{\displaystyle \sum_{i = 1}^m }\left\Vert\nabla\varphi_i\right\Vert_{L^2(\Omega)} -\beta_0.
\]

Thus,
\[
D{\displaystyle \sum_{i = 1}^m }\lambda_{0i}\geqslant \left({\displaystyle \sum_{i = 1}^m }\left\Vert\nabla\varphi_i\right\Vert^2_{L^2(\Omega)}\right)^{\frac{1}{2}}\left(c_0\left({\displaystyle \sum_{i = 1}^m }\left\Vert\nabla\varphi_i\right\Vert^2_{L^2(\Omega)}\right)^{\frac{1}{2}}  - M\sqrt{m}\right) -\beta_0.
\]

According to (\ref{4.2}), by the variational characterization of $\lambda_{1}(-\Delta,\Omega)$, we have that
\[
{\displaystyle \sum_{i = 1}^m }\left\Vert\nabla\varphi_i\right\Vert^2_{L^2(\Omega)}\geqslant \lambda_{1}(-\Delta,\Omega){\displaystyle \sum_{i = 1}^m }\left\Vert\varphi_i\right\Vert^2_{L^2(\Omega)}=\lambda_{1}(-\Delta,\Omega).
\]

Furthermore, by Faber-Krahn's inequality (see page 280 of \cite{Lo}), it is well known that
\[
\lambda_{1}(-\Delta,\Omega)\geqslant \Sigma\vert B_1\vert^{\frac{2}{n}}\vert\Omega\vert^{-\frac{2}{n}}.
\]

Therefore, for $\vert\Omega\vert$ verifying (\ref{4.21}), we obtain (\ref{4.1}) holds. 

Now, we take $\sigma_0=(1,1,\ldots,1)\in (0,\infty)^{m-1}$. Then 
\[
\lambda_*=\left[\theta_*(\sigma_0)\right]^{\displaystyle{\sum_{j=1}^{m}\prod_{i=1}^j\alpha_i}}.
\]
Therefore, if in addition $\Omega$ satisfies (\ref{01}), since $\theta_*(\sigma_0) \theta_*(\sigma_0)^{\alpha_1} \ldots \theta_*(\sigma_0)^{\alpha_1\ldots\alpha_{m-1}}=\lambda_*$, we have 
\[
\lambda_*\geqslant \theta_*(\sigma_0)+(\theta_*(\sigma_0))^{\alpha_1}+\ldots+(\theta_*(\sigma_0))^{\alpha_1\ldots\alpha_{m-1}}> m.
\]
Finally, from (\ref{4.1}), we obtain inequality (\ref{02}). This ends the proof of Theorem \ref{teo1}. \ \rule {1.5mm}{1.5mm}\\

{\bf Proof of Corollary \ref{EQ1}} Note that, since $\Omega$ satisfies \eqref{hip5} and \eqref{hip6}, we have the operator ${\cal L}_i$ satisfy {\bf (SMP)} in $\Omega$, for $i = 1,\ldots,m$.
By characterization (\ref{carac}), we have (ii) implies $\Lambda \in \overline{\mathcal{R}_1} \setminus {\bf \Lambda_1}$, that is, $\Lambda\geqslant 0$. Moreover, the equivalence between (ii) and (iii) also follow of (\ref{carac}). Then, it suffices to show that (i) implies (ii).

Again by (\ref{carac}) the conclusion is immediate in the cases that $\lambda_j = 0$ for some $j\in\{1,\ldots,m\}$. So, it suffices to consider $\Lambda \in (0, \infty)^m$. In this case, we take $\Lambda_0=(\theta_*(\sigma), \theta_*(\sigma) \sigma)\in {\bf \Lambda_1}$, with $$\sigma=\left(\frac{\lambda_2}{\lambda_1},\ldots,\frac{\lambda_{m}}{\lambda_1}\right)\ \ \in\ \ (0, \infty)^{m-1}.$$ Since $\Omega$ satisfies (\ref{4.3}), by Theorem \ref{teo1}
\[
\theta_*(\sigma)\left(1+\frac{\lambda_2}{\lambda_1}+\ldots+\frac{\lambda_{m}}{\lambda_1}\right)>{\displaystyle \sum_{i = 1}^m }\lambda_{i}.
\]
Therefore, $\theta_*(\sigma)>\lambda_1$ and thus $\Lambda_0>\Lambda$. Then, $\Lambda \in \overline{{\cal R}_1} \setminus {\bf \Lambda_1}$. Applying again the characterization (\ref{carac}), we have {\bf (WMP)} associated to \eqref{1.3} holds in $\Omega$. This concludes the proof of Corollary \ref{EQ1}. \ \rule {1.5mm}{1.5mm}

\section{Proof of Theorem \ref{LE}}

For the proof of the Theorem \ref{LE} it will be necessary the use of the following inequalities.

$\bullet$ For $n = 1$ there is a constant $C_1 > 0$ such that,

\begin{equation}\label{sup}
||u||_{L^\infty(\Omega)} \leqslant C_1 ||\nabla u||_{L^2(\Omega)}\, ,
\end{equation}
for any $u \in H_0^{1}(\Omega)$. 

$\bullet$ For $n = 2$, we invoke the Moser-Trudinger inequality, which states (see \cite{Moser} and \cite{Tr}) that

\begin{equation}\label{MT}
\int_\Omega e^{4 \pi u^2} dx \leqslant \xi |\Omega|
\end{equation}
for all $u \in H_0^{1}(\Omega)$ such that $||\nabla u||_{L^2(\Omega)} \leqslant 1$, where the positive constant $\xi$ independent of $\Omega$. Then
\begin{equation} \label{3}
||u||^2_{L^p(\Omega)} \leqslant C_2 |\Omega|^{2/p} ||\nabla u||^2_{L^2(\Omega)}\, ,
\end{equation}
because $|u|^p \leqslant C(p) e^{4 \pi u^2}$ for any $p \geqslant 1$ and instead $u$, we use $u/||\nabla u||_{L^2(\Omega)}$ in the inequality (\ref{MT}), with $C_2 = (C(p) \xi)^{2/p}$ for any $p \geqslant 1$.

$\bullet$ For $n \geqslant 3$, we invoke the sharp Sobolev embedding for any $u \in H_0^{1}(\Omega)$,

\begin{equation} \label{4}
||u||^2_{L^{\frac{2n}{n-2}}(\Omega)} \leqslant C_n ||\nabla u||^2_{L^2(\Omega)},
\end{equation}
where an explicit expression of $C_n>0$ depending only on $n$ was proved in \cite{Au, Ta}.

We first focus on the claim (iii). Let $\varphi$ a positive eigenfunction corresponding to $\Lambda_0\in {\bf \Lambda_1}$. So, multiplying by $\varphi_1$ the problem
\[
-{\cal L}_1\varphi_1 = \lambda_{01}\rho_1(x)\varphi_2^{\alpha_1}
\]
and integrating in $\Omega$, we have

\[
\lambda_{01}\int_{\Omega}\rho_1(x)\varphi_2^{\alpha_1}\varphi_1\geqslant c_0\left\Vert\nabla\varphi_1\right\Vert_{L^2(\Omega)}^2-\int_\Omega\langle b^1,\nabla\varphi_1\rangle\varphi_1-\inf_\Omega c_1\int_\Omega\varphi_1^2.
\]
Note that, by applying H\"{o}lder and Young's inequalities and (\ref{hip0}), we get
\begin{eqnarray*}
\left\vert\int_\Omega\langle b^1,\nabla\varphi_1\rangle\varphi_1\right\vert &\leqslant &\frac{c_0}{2}\left\Vert\nabla\varphi_1\right\Vert_{L^{2}(\Omega)}^2+\vert\Omega\vert^{\frac{p(\alpha_1-1)+2}{2p(\alpha_1+1)}}\Vert b^1\Vert_\infty\left\Vert\varphi_1\right\Vert_{L^{\gamma+1}(\Omega)}^2
\end{eqnarray*}
and
\[
\int_{\Omega}\rho_1(x)\varphi_2^{\alpha_1}\varphi_1\leqslant C\left(\Vert \varphi_1\Vert_{L^{\gamma+1}(\Omega)}^2+\Vert \varphi_2\Vert_{L^{\gamma+1}(\Omega)}^{2\alpha_1}\right),
\]
where $\gamma=\frac{p\alpha_1+1}{p-1}$. Note also that
\[
-\inf_\Omega c_1\int_\Omega\varphi_1^2\geqslant -\vert\Omega\vert^{\frac{p(\alpha_1-1)+2}{2p(\alpha_1+1)}}\vert\inf_\Omega c_1\vert\Vert\varphi_1\Vert_{L^{\gamma+1}(\Omega)}^2
\]
from the H\"{o}lder's inequality. Thus,

\begin{equation}\label{ine1}
\lambda_{01}C\left(\Vert \varphi_1\Vert_{L^{\gamma+1}(\Omega)}^2+\Vert \varphi_2\Vert_{L^{\gamma+1}(\Omega)}^{2\alpha_1} \right)\geqslant\frac{c_0}{2}\left\Vert\nabla\varphi_1\right\Vert_{L^2(\Omega)}^2-c\vert\Omega\vert^{\frac{p(\alpha_1-1)+2}{2p(\alpha_1+1)}}\left\Vert\varphi_1\right\Vert_{L^{\gamma+1}(\Omega)}^2.
\end{equation}
Similarly, for $i=2,\ldots,m$ it follows from the problem $-{\cal L}_i\varphi_i = \lambda_{0i}\rho_i(x)\varphi_{i+1}^{\alpha_i}$ that
\[
\lambda_{0i}C\vert\Omega\vert^{\frac{(\alpha_1-\alpha_i)(n-1)}{n(\alpha_1+1)}}\left(\Vert \varphi_{i+1}\Vert_{L^{\gamma+1}(\Omega)}^{2\alpha_i}+\Vert \varphi_i\Vert_{L^{\gamma+1}(\Omega)}^{2} \right)\geqslant\frac{c_0}{2}\left\Vert\nabla\varphi_i\right\Vert_{L^2(\Omega)}^2-c\vert\Omega\vert^{\frac{p(\alpha_1-1)+2}{2p(\alpha_1+1)}}\left\Vert\varphi_i\right\Vert_{L^{\gamma+1}(\Omega)}^2.
\]
Then, from the (\ref{hip1}), we derive

\begin{equation}\label{ine2}
\lambda_{0i}^{\Pi_{j = 1}^{i-1} \alpha_j} C\left(\Vert \varphi_{i+1}\Vert_{L^{\gamma+1}(\Omega)}^{2\Pi_{j = 1}^{i} \alpha_j}+\Vert \varphi_i\Vert_{L^{\gamma+1}(\Omega)}^{2\Pi_{j = 1}^{i-1} \alpha_j} \right)\geqslant\tilde{c_0}\left\Vert\nabla\varphi_i\right\Vert_{L^2(\Omega)}^{2\Pi_{j = 1}^{i-1} \alpha_j}-\overline{c_0}\vert\Omega\vert^{\frac{p\beta(\alpha_1-1)+2\beta}{2p(\alpha_1+1)}}\left\Vert\varphi_i\right\Vert_{L^{\gamma+1}(\Omega)}^{2\Pi_{j = 1}^{i-1} \alpha_j}.
\end{equation}
Therefore, adding up inequalities (\ref{ine1}) and (\ref{ine2}) for $i = 2,\ldots,m$, we have
\[
2C\max\{\lambda_{01},\lambda_{02}^{\alpha_{1}},\ldots,\lambda_{0m}^{\alpha_1\ldots\alpha_{m-1}}\}\geqslant \tilde{c_0}\left(\frac{{\displaystyle \sum_{i = 1}^m }\left\Vert\nabla\varphi_i\right\Vert_{L^2(\Omega)}^{\frac{2}{\Pi_{j = i}^m \alpha_j}}}{{\displaystyle \sum_{i = 1}^m }\left\Vert\varphi_i\right\Vert_{L^{\gamma+1}(\Omega)}^{\frac{2}{\Pi_{j = i}^m \alpha_j}}}\right)- \overline{c_0}\vert\Omega\vert^{\frac{p\beta(\alpha_1-1)+2\beta}{2p(\alpha_1+1)}}.
\]

For $i=1,\ldots,m$, by the sharp Sobolev embedding (\ref{4}), the variational characterization of $\lambda_1(-\Delta,\Omega)$, interpolation inequality and using again (\ref{hip1}), we have
\[
\frac{\Vert\nabla\varphi_i\Vert_{L^2(\Omega)}^{\frac{2}{\Pi_{j = i}^m \alpha_j}}}{\Vert\varphi_i\Vert_{L^{\gamma+1}(\Omega)}^{\frac{2}{\Pi_{j = i}^m \alpha_j}}}\geqslant\left(\frac{1}{C_n^{1-\theta}}\lambda_1(-\Delta,\Omega)^\theta\right)^{\frac{1}{\Pi_{j = i}^m \alpha_j}}\geqslant\frac{1}{C_n^{\beta-\theta\beta}}\lambda_1(-\Delta,\Omega)^{\theta\beta},
\]
where
\[
\frac{1}{\gamma+1}=\frac{\theta}{2}+\frac{(1-\theta)(n-2)}{2n}.
\]

Then,

\[
\frac{{\displaystyle \sum_{i = 1}^m }\left\Vert\nabla\varphi_i\right\Vert_{L^2(\Omega)}^{\frac{2}{\Pi_{j = i}^m \alpha_j}}}{{\displaystyle \sum_{i = 1}^m }\left\Vert\varphi_i\right\Vert_{L^{\gamma+1}(\Omega)}^{\frac{2}{\Pi_{j = i}^m \alpha_j}}}\geqslant\frac{1}{C_n^{\beta-\theta\beta}}\lambda_1(-\Delta,\Omega)^{\theta\beta}.
\]
So, by Faber-Krahn's inequality, we obtain
\[
\max\{\lambda_{01},\lambda_{02}^{\alpha_{1}},\ldots,\lambda_{0m}^{\alpha_1\ldots\alpha_{m-1}}\}\geqslant\frac{1}{2C}\left(\frac{1}{C_n^{\beta-\theta\beta}}\tilde{c_0}\Sigma^{\theta\beta}\vert B_1\vert^{\frac{2\theta\beta}{n}}\vert\Omega\vert^{-\frac{2\theta\beta}{n}}-\overline{c_0}\vert\Omega\vert^{\frac{p\beta(\alpha_1-1)+2\beta}{2p(\alpha_1+1)}}\right).
\] 

Now, let $\sigma_0=(1,1,\ldots,1)\in (0,\infty)^{m-1}$. Then 
\[
\lambda_*=\left[\theta_*(\sigma_0)\right]^{\displaystyle{\sum_{j=1}^{m}\prod_{i=1}^j\alpha_i}}.
\]

Thus, using that $\Omega$ satisfies (\ref{h1}) and the fact that $\theta_*(\sigma_0) \theta_*(\sigma_0)^{\alpha_1} \ldots \theta_*(\sigma_0)^{\alpha_1\ldots\alpha_{m-1}}=\lambda_*$, we obtain
\[
\lambda_*\geqslant \max\{\theta_*(\sigma_0),\theta_*(\sigma_0)^{\alpha_{1}},\ldots,\theta_*(\sigma_0)^{\alpha_1\ldots\alpha_{m-1}}\}> 1.
\]

Then (\ref{8.108}) occurs, and so we conclude the proof of the assertion (iii). Finally, the proof of claims (i) and (ii) follow in a similar way, instead of inequality (\ref{4}), we invoke respectively the inequalities (\ref{sup}) and (\ref{3}). Furthermore, for the claim (ii), instead of hypotheses (\ref{h1}), we use (\ref{h2}). This finishes the proof.  \ \rule {1.5mm}{1.5mm}\\

\section{Proof of Theorem \ref{upper}}

Let $B_r=\{x \in \R^n: |x| < r\}$ with $r=R/2$. We consider the function $\tau$ given by

\[
\tau(x) =\frac{1}{4}(r^2-\vert x\vert^2)^2.
\]
Note that $\tau\in C^2(\overline{B}_r)$ with $0\leqslant\tau\leqslant\frac{r^4}{4}$ in $\overline{B}_r$, $\tau=0$ on $\partial B_r$, $\tau(0)=\frac{r^4}{4}$ and $\tau>0$ in $B_r$. We assertion that for $i=1,\ldots,m$
\begin{equation}\label{1.17}
\sup_{B_r}\left(-\frac{{\cal L}_i \tau}{\rho_i(x) \tau^{\alpha_i}}\right)\leqslant\frac{A_i}{R^{4\alpha_i-2}},
\end{equation}
where $A_i$ is given in (\ref{hip3}). Indeed, for each $i=1,\ldots,m$, a direct calculation ensures that
\begin{eqnarray*}
-\frac{{\cal L}_i \tau}{4^{\alpha_i} \tau^{\alpha_i}} &=& \frac{1}{4^{\alpha_i} \tau^{\alpha_i}} \left[ (r^2 - |x|^2) \sum_{\kappa =1}^n a^i_{\kappa\kappa} - 2 \sum_{\kappa,l =1}^n a^i_{\kappa l} x_\kappa x_l + \sum_{\kappa =1}^{n} b^i_\kappa x_\kappa (r^2 - |x|^2) - c_i \tau  \right] \nonumber\\
&\leqslant &\frac{n C_0+b_0 r}{(r^2-\vert x\vert^2)^{2\alpha_i - 1}}-\frac{2c_0\vert x\vert^2}{(r^2-\vert x\vert^2)^{2\alpha_i}}+\frac{b_0}{4(r^2-\vert x\vert^2)^{2\alpha_i - 2}} \nonumber \\
&\leqslant &\frac{r^2(4n C_0+b_0(4r + r^2))-\vert x\vert^2(4n C_0+b_0(4r + r^2)+2c_0)}{(r^2-\vert x\vert^2)^{2\alpha_i}}.
\end{eqnarray*}
Note that, in the region of $B_r$ where
\[
|x|^2 \left(4nC_0 + b_0(4r + r^2) + 2c_0 \right) > r^2 \left( 4n C_0 + b_0(4r+r^2) \right),
\]
we obtain
\[
-\frac{{\cal L}_i \tau}{\rho_i(x) \tau^{\alpha_i}} \leqslant 0.
\]
In the remainder of $B_r$, since $R < 1$, we have 
\begin{eqnarray*}
-\frac{{\cal L}_i \tau}{4^{\alpha_i} \tau^{\alpha_i}} &\leqslant& \frac{4n C_0 + b_0(4r + r^2)}{(r^2-\vert x\vert^2)^{2\alpha_i - 1}} \\
&\leqslant& \frac{1}{\left(2c_0r^2\right)^{2\alpha_i - 1}} \left( 4n C_0+b_0(4r + r^2) \right) \left(4n C_0+b_0(4r + r^2) + 2c_0 \right)^{2\alpha_i - 1}\\
&\leqslant& \frac{1}{\left(2c_0r^2\right)^{2\alpha_i - 1}} \left( 4n C_0+\frac{9}{4} b_0 \right) \left(4n C_0+\frac{9}{4}b_0 + 2c_0 \right)^{2\alpha_i - 1}
\end{eqnarray*}
for $\alpha_i \geqslant 1/2$ and 
\[
-\frac{{\cal L}_i \tau}{4^{\alpha_i} \tau^{\alpha_i}} \leqslant \frac{4n C_0 + b_0(4r + r^2)}{(r^2-\vert x\vert^2)^{2\alpha_i - 1}} \leqslant \frac{4n C_0 + b_0(4r + r^2)}{(r^2)^{2\alpha_i - 1}}\leqslant \frac{4n C_0 + \frac{9}{4}b_0}{(r^2)^{2\alpha_i - 1}}
\]
for $\alpha_i < 1/2$. Therefore,
\[
-\frac{{\cal L}_i \tau}{\rho_i(x) \tau^{\alpha_i}} \leqslant \frac{4^{\alpha_i}}{\varepsilon_0} \left( -\frac{{\cal L}_i \tau}{4^{\alpha_i} \tau^{\alpha_i}} \right) \leqslant \frac{A_i}{R^{4\alpha_i - 2}},
\]
and thus \eqref{1.17} holds. Set
\[
\lambda_i:=\sup_{B_r}\left(-\frac{{\cal L}_i \tau}{\rho_i(x) \tau^{\alpha_i}}\right).
\]
By inequality (\ref{hip2}), we obtain

\[
{\cal L}_i \tau(0) = - r^2 \sum_{\kappa =1}^n a^i_{\kappa\kappa}(0) + \frac{c_i(0)}{4} r^4 \leqslant r^2 (- nc_0 + \frac{c_i(0)}{4} r^2) < 0
\]
for $i =1,\ldots, m$. Note that, from hypotheses of $R$, we get $\lambda_i > 0$.

We assertion that $\Lambda=(\lambda_1,\ldots,\lambda_m)\notin  \overline{{\cal R}_1(\Omega)}\setminus{\bf \Lambda_1}(\Omega)$. Assume by contradiction that $\Lambda \in  \overline{{\cal R}_1(\Omega)}\setminus{\bf \Lambda_1}(\Omega)$. Note that for $i =1,\ldots, m$, we have
\[
\begin{array}{llll}
-{\cal L}_i (-\tau) - \lambda_i\rho_i(x) |-\tau|^{\alpha_i-1}(-\tau) = -{\cal L}_i (-\tau) - \sup\limits_{B_r}\left(-\frac{{\cal L}_i\tau}{\rho_i \tau^{\alpha_i}}\right)\rho_i(x)|-\tau|^{\alpha_i-1}(-\tau) \geqslant 0
\end{array}
\ {\rm in}\ B_r\ {\rm a.e.}
\]
and $-\tau = 0$ in $\partial B_r$. Then, by Theorem 6.1 of \cite{LM1}, we derive $\Lambda \in  \overline{{\cal R}_1(B_r)}\setminus{\bf \Lambda_1}(B_r)$ because $B_r \subset \Omega$.  Invoking {\bf (WMP)} in $B_r$ we are led to the contradiction $\tau \leqslant 0$ in $B_r$. This conclude the proof of claim. Thus, $\Lambda_0=(\theta_*(\sigma_0), \theta_*(\sigma_0) \sigma_0)\leqslant\Lambda$, where $$\sigma_0=\left(\frac{\lambda_2}{\lambda_1},\ldots,\frac{\lambda_{m}}{\lambda_1}\right)\ \ \in\ \ (0, \infty)^{m-1}.$$

Therefore, from \eqref{1.17}, we obtain
\begin{eqnarray*}
\lambda_* &\leqslant &\frac{A_1}{R^{4\alpha_1-2}}\left(\frac{A_2}{R^{4\alpha_2-2}}\right)^{\alpha_1} \ldots \left(\frac{A_m}{R^{4\alpha_m-2}}\right)^{\alpha_1\ldots\alpha_{m-1}}\\
&= & A_1 A^{\alpha_1}_2 \ldots A^{\alpha_1\ldots\alpha_{m-1}}_m R^{h}=ER^{ h},
\end{eqnarray*}
where $E=E(n, \varepsilon_0, \alpha, c_0, C_0, b_0)>0$ and $h=-2\displaystyle{\sum_{j=1}^{m}\prod_{i=1}^j\alpha_i}$. This ends the proof.\ \rule {1.5mm}{1.5mm}\\

\section{Proof of Theorem \ref{teo2}}

The natural space where we look for solutions to problem \eqref{poly} (see \cite{GGS}) is defined as follows:
$$
H^m_{\vartheta}(\Omega):=\left\{v\in W^{m,2}(\Omega)\colon\ \Delta^jv =0\ \text{on}\ \partial\Omega\ \text{for}\ j<\frac{m}{2}\right\},
$$
that inherits the scalar product and the Hilbert space structure from $W^{m,2}(\Omega),$ then
 $$
(v,\phi)_{H^m_{\vartheta}}=\sum_{|\zeta|\leqslant m}(D^\zeta v,D^\zeta \phi)_{L^2(\Omega)},\ \zeta=(\zeta_1,\cdots,\zeta_n), \quad\text{and}\quad \|v\|_{H^m_{\vartheta}}=(v,v)^{1/2}_{H^m_\vartheta}.
$$

Consider the principal $(m-1)$-hypersurface ${\bf \Lambda_1}={\bf \Lambda_1}(-\Delta,\Omega,\rho)$ corresponding to following GLE system

\begin{equation}\label{pol}
\left\{
\begin{array}{rlllr}
-\Delta u_1 &=& \lambda_{1} \rho_1(x) u_2 & {\rm in} & \Omega, \\
-\Delta u_2 &=& \lambda_2 u_3 & {\rm in} & \Omega, \\
& \vdots & & \vdots &\\
-\Delta u_m &=& \lambda_m u_1 & {\rm in} & \Omega, \\
u&=&0 & {\rm on} & \partial \Omega,
\end{array}\right.
\end{equation}
where $\Delta$ is the classical Laplace operator and $\rho=(\rho_1,1,\ldots,1) \in L^p(\Omega;\R^m)$, with $p>n$. We take $\Lambda_0=(\lambda_{01},\ldots,\lambda_{0m})\in {\bf \Lambda_1}$, with $\lambda_{01}=\lambda_*$ and $\lambda_{0i}=1$, $i = 2,\ldots, m$. Then, for all principal eigenfunction $\varphi=(\varphi_1,\ldots,\varphi_m)$ of (\ref{pol}) corresponding to $\Lambda_0$, we have $\varphi_2$ is a solution of the problem (\ref{poly}). Furthermore, $\lambda_{01}$ is a positive eigenvalue of the problem (\ref{poly}) that admits eigenfunction $v\in W^{2m,p}(\Omega)$ satisfying $v, (-\Delta)^j v>0$ in $\Omega$, $j = 1,\ldots, m-1$. 
We take $\lambda_1((-\Delta)^{m},\Omega,\rho_1):=\lambda_{01}=\lambda_*$.

Suppose that there exists $0<\lambda<\lambda_1((-\Delta)^{m},\Omega,\rho_1)$ eigenvalue to \eqref{poly}, that is, there is a nontrivial function $z\in H^m_{\vartheta}(\Omega)$ satisfying

  \[
    (z,\phi)_{H^m_\vartheta}=\lambda\int_{\Omega}\rho_1(x)z\phi\,dx\quad \forall \phi\in H^m_{\vartheta}(\Omega)
  \]
  and $z=(-\Delta)^j z=0$ on $\partial\Omega$, $j = 1,\ldots, m-1$.

Using Sobolev embedding, Theorem 2.20 of \cite{GGS} and a bootstrap argument, we obtain $z\in W^{2m,p}(\Omega)$.

Now, we take $u_2=z,\ u_3=-\Delta u_2,\ u_4=-\Delta u_3=(-\Delta)^2 u_2,\ldots, u_m=(-\Delta)^{m-2} u_2$ and $u_1=(-\Delta)^{m-1} u_2$. Then, $u=(u_1,\ldots,u_m)\in W^{2,p}(\Omega;\R^m)$ is a strong solution of the problem
 
\[
\left\{
\begin{array}{rlllr}
-\Delta u_1 &=& \lambda \rho_1(x) u_2 & {\rm in} & \Omega, \\
-\Delta u_2 &=& u_3 & {\rm in} & \Omega, \\
& \vdots & & \vdots &\\
-\Delta u_{m-1} &=& u_m & {\rm in} & \Omega, \\
-\Delta u_m &=& u_1 & {\rm in} & \Omega, \\
u&=&0 & {\rm on} & \partial \Omega.
\end{array}\right.
\]

Thus $(\lambda,1,\cdots,1)\in {\cal R}_1$ (see Figure 2 for $m=2$). 
\begin{figure}[ht]
\centering
\includegraphics[scale=1]{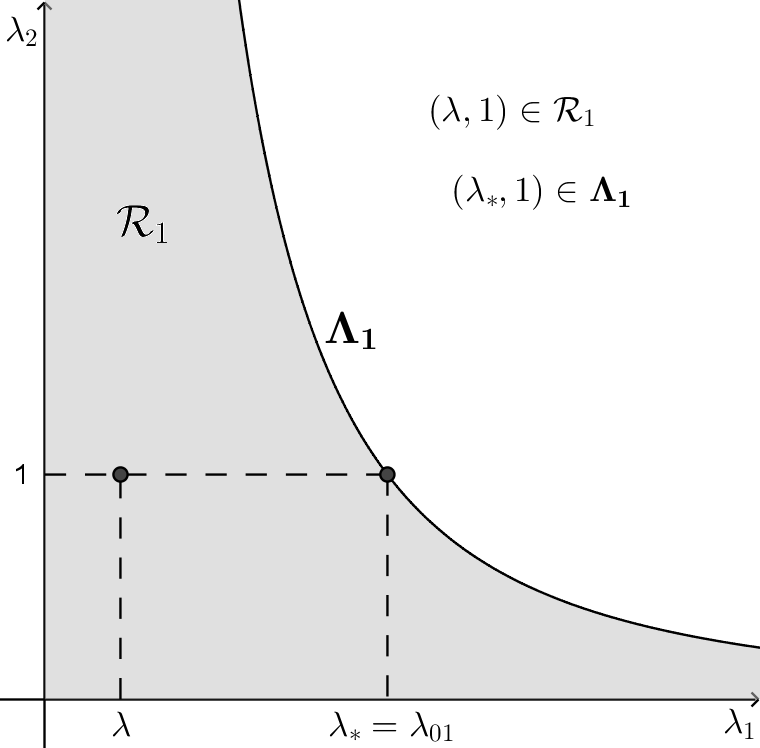}
\caption{For $m=2$, ${\bf \Lambda_1}=\left\{(\lambda_1,\lambda_2) \in (0, \infty)^2 :\, \lambda_1 \lambda_2 = \lambda_*\right\}$ and ${\cal R}_1 := \{t \Lambda: 0 < t < 1\ {\rm and}\ \Lambda \in {\bf \Lambda_1}\}$.}
\end{figure}

But it's a contradiction, because by Theorem 6.1 of \cite{LM1}, there isn't any eigenvalue to (\ref{1.3}) in the region $(\overline{{\cal R}}_1 \setminus {\bf \Lambda_1}) \cup - (\overline{{\cal R}}_1 \setminus {\bf \Lambda_1})$.

From the simplicity of ${\bf \Lambda_1}$, it follows that $\lambda_1((-\Delta)^{m},\Omega,\rho_1)$ is simple. \ \rule {1.5mm}{1.5mm}\\

\section{Bounds for eigenvalues of (\ref{poly}) and proof of Theorems \ref{teo4} and \ref{cor3}}

In this section, we consider the principal eigenvalue $\lambda_1((-\Delta)^m,\Omega,\rho_1)$ of the problem (\ref{poly}) obtained through the GLE system (\ref{pol}). We separate the estimates for $\lambda_1((-\Delta)^m,\Omega,\rho_1)$ into three items below.\\

{\bf 1) Lower estimate for $\lambda_1((-\Delta)^m,\Omega,\rho_1)$ when $\rho_1 \in L^p(\Omega)$, with $p>n$}\\

Let $C=\max\{1,\Vert\rho_1\Vert_{L^p(\Omega)},\vert\Omega\vert^{\frac{1}{p}}\}$ and $\Omega$ satisfying (\ref{hip1}).

$\bullet$ For $n=1$ and $\Omega$ satisfying
\[
\vert\Omega\vert< \frac{\vert B_1\vert\Sigma^{\frac{1}{2}}}{(C(C_1)^{1/p})^\frac{1}{2\theta}},
\]
by Theorem \ref{LE}, we get
\[
\lambda_1((-\Delta)^m,\Omega,\rho_1)\geqslant\frac{1}{4C(C_1)^{1/p}}\left(\Sigma\vert B_1\vert^{2}\vert\Omega\vert^{-2}\right)^{1-1/p},
\]
where $\theta=\frac{p-1}{p}$.

$\bullet$ For $n=2$ and $\Omega$ satisfying
\[
\vert\Omega\vert< \frac{\left(\vert B_1\vert^{2}\Sigma^{2}\right)^\theta}{(C(C_2)^{2/3})^2},
\]
by Theorem \ref{LE}, we derive
\[
\lambda_1((-\Delta)^m,\Omega,\rho_1)\geqslant\frac{1}{4C(C_2)^{2/3}}\left(\Sigma^{\frac{1}{3}}\vert B_1\vert^{\frac{1}{3}}\vert\Omega\vert^{-\frac{1}{2}}\right),
\]
where $\theta=\frac{1}{3}$.

$\bullet$ For $n\geqslant 3$ and $\Omega$ satisfying
\[
\vert\Omega\vert< \frac{\vert B_1\vert\Sigma^{\frac{n}{2}}}{(C(C_n)^{1-\theta})^\frac{n}{2\theta}},
\]
by Theorem \ref{LE}, we obtain
\[
\lambda_1((-\Delta)^m,\Omega,\rho_1)\geqslant\frac{1}{4C(C_n)^{1-\theta}}\left(\Sigma\vert B_1\vert^{\frac{2}{n}}\vert\Omega\vert^{-\frac{2}{n}}\right)^{\theta},
\]
where $\frac{p-1}{2p}=\frac{\theta}{2}+\frac{(1-\theta)(n-2)}{2n}$.\\

{\bf 2) Lower estimate for $\lambda_1((-\Delta)^m,\Omega,\rho_1)$ when $\rho_1 \in L^\infty(\Omega)$}\\

In this case, by Theorem \ref{teo1}, for any $\Omega$ satisfying 
\begin{equation}\label{h01}
\vert\Omega\vert\leqslant \frac{\vert B_1\vert\Sigma^{\frac{n}{2}}}{(mD)^\frac{n}{2}},
\end{equation}
we have
\begin{equation}\label{i1}
\lambda_1((-\Delta)^m,\Omega,\rho_1)\geqslant\frac{1}{D}\left(\Sigma\vert B_1\vert^{\frac{2}{n}}\vert\Omega\vert^{-\frac{2}{n}}\right),
\end{equation}
where $D=\max\{1,\Vert\rho_1\Vert_{L^\infty(\Omega)}\}$.\\

{\bf 3) Upper estimate for $\lambda_1((-\Delta)^m,\Omega,\rho_1)$ when $\rho_1 \in L^\infty(\Omega)$}\\

By Theorem \ref{upper}, if $\rho_1 \in L^\infty(\Omega)$, with $\rho_1 \geqslant \varepsilon_0$ in $\Omega$ for some constant $\varepsilon_0 > 0$, $0 \in \Omega$ and $0<R < 1$ such that $\overline{B}_R \subset \Omega$, we get
\begin{equation}\label{i2}
\lambda_1((-\Delta)^m,\Omega,\rho_1)\leqslant \left(\dfrac{64n(2n+1)}{\varepsilon_0R^2}\right)^m.
\end{equation}

\begin{rem}
Under the hypotheses of Theorem \ref{teo4}, we claim that given $\mu\geqslant m$, there are $s,R,\varepsilon_0>0$ such that $\lambda_1((-\Delta)^m,B_s,\rho_1)=\mu$. In fact,
we take $0<R<1$, $R<s$ and $\rho_1 \geqslant \varepsilon_0$ in $B_s$, with $$\varepsilon_0=\dfrac{64n(2n+1)}{\mu^{\frac{1}{m}}R^2}.$$ We also take $s$ satisfying
\[
s^2=\frac{\Sigma}{\mu D}.
\]
Then $B_s$ satisfies (\ref{h01}). Therefore, by inequalities in (\ref{i4}), we obtain $\lambda_1((-\Delta)^m,B_s,\rho_1)=\mu$.
\end{rem}

{\bf Proof of Theorem \ref{teo4}} The proof of this result follows from inequalities (\ref{i1}) and (\ref{i2}) with $\Omega=B_s$.  \ \rule {1.5mm}{1.5mm}\\

{\bf Proof of Theorem \ref{cor3}} If $\lambda=0$ then the result follows from Theorem 1.1 of \cite{LM1}. Suppose $\lambda>0$ and take $u_2=v,\ u_3=-\Delta u_2,\ u_4=-\Delta u_3=(-\Delta)^2 u_2,\ldots, u_m=(-\Delta)^{m-2} u_2$ and $u_1=(-\Delta)^{m-1} u_2$. Then $u=(u_1,\ldots,u_m)\in W^{2,p}(\Omega;\R^m)$ is a strong solution of the problem

\[
\left\{
\begin{array}{rlllr}
-\Delta u_1 &\geqslant & \lambda \rho_1(x) u_2 & {\rm in} & \Omega, \\
-\Delta u_2 &=& u_3 & {\rm in} & \Omega, \\
& \vdots & & \vdots &\\
-\Delta u_m &=& u_1 & {\rm in} & \Omega, \\
u&=&0 & {\rm on} & \partial \Omega.
\end{array}\right.
\]

Since $\Lambda=(\lambda,1,\ldots,1)>0$, by Corollary \ref{cor2}, we obtain $u>0$ in $\Omega$ whenever $\vert\Omega\vert<\eta_0$, i.e., $v,(-\Delta)^j v>0$ in $\Omega$, $j = 1,\ldots, m-1$ provided that $\vert\Omega\vert<\eta_0$. This finish the proof of Theorem \ref{cor3}.\ \rule {1.5mm}{1.5mm}\\

{\bf Acknowledgments:} S. Bahrouni has been supported by FAPESP Proc 2023/04515-7. E.J.F. Leite has been supported by CNPq/Brazil (PQ 316526/2021-5). G.F. Madeira has been supported by FAPESP Proc 2022/16407-1.\\

{\bf Data availability}
Data sharing not applicable to this article as no datasets were generated or analyzed during the current study.\\

%{\bf Conflict of interest}
%The authors declare that they have no conflict of interest.\\
%
%{\bf Ethics approval}
%The research does not involve humans and/or animals. The authors declare that there are no
%ethics issues to be approved or disclosed.\\

\end{document}